\documentclass[11pt]{article}
\usepackage{amsmath}
\usepackage{amssymb}
\usepackage{epsf}
\usepackage{amscd}
\usepackage[all]{xy}
\usepackage[francais,english]{babel}

\input xy
\xyoption{all}

\newtheorem{theorem}{Theorem}[section]
\newtheorem{propo}[theorem]{Proposition}
\newtheorem{coro}[theorem]{Corollary}
\newtheorem{lemme}[theorem]{Lemma}

\newenvironment{theo}{\begin{theorem} \sf}{\end{theorem}}
\newenvironment{prop}{\begin{propo} \sf}{\end{propo}}
\newenvironment{cor}{\begin{coro} \sf}{\end{coro}}
\newenvironment{lemma}{\begin{lemme} \sf}{\end{lemme}}

\def\ot{\otimes}

\advance\oddsidemargin by- 1.4cm
\numberwithin{equation}{section}

\title{\textsc{Cohomology theories of Hopf bimodules and cup-product.}}

\author{Rachel Taillefer\thanks{Laboratoire G.T.A., UPRES A 5030,
D\'epartement de Math\'ematiques CC 051, Universit\'e Montpellier II,
34095 Montpellier Cedex 5, France. 
E-mail: taillefr\at math.univ-montp2.fr}}

\date{}

\begin{document}

\maketitle


\selectlanguage{english}
\begin{abstract} 
Given a Hopf algebra $A,$ there exist various cohomology theories for
the category of Hopf bimodules over $A,$ introduced by M. Gerstenhaber
and S.D. Schack, and by C. Ospel. We prove, when $A$ is finite dimensional, that they are all equal to the $Ext$ functor on the module category of an associative algebra associated to $A,$ described by C. Cibils and M. Rosso. We also give an expression for a cup-product in the cohomology defined by C. Ospel, and prove that it corresponds to the Yoneda product of extensions.
\end{abstract}


\section{Introduction}

Let $A$ be a finite dimensional Hopf algebra over a field $k,$ and let $M$ and $N$ be  Hopf bimodules over $A.$ We shall consider various cohomology theories associated to these objects. In the first place, M. Gerstenhaber and S.D. Schack have defined, in~\cite{GS1}, a cohomology $H^*_b(A,A)$ for $A.$ It arose when studying the deformations of $A$ (however, this cohomology does not characterize the deformations of $A$ in the category of Hopf algebras, but in a larger class of objects which includes Hopf algebras and Drinfel'd algebras). Next, in~\cite{GS2}, M. Gerstenhaber and S.D. Schack extended this cohomology to a cohomology theory $H_{GS}^*(M,N)$ of $M$ into $N.$ Lastly, C. Ospel defined in his thesis~\cite{O} a cohomology $H_{A4}^* (M,A)$ of $M$ with coefficients in $A,$ introducing a double complex which involves both  Hochschild cohomology complexes and  Cartier cohomology complexes~(\cite{C}).

In this paper, we will show that these cohomologies are all the
same. To this end, we will introduce an associative algebra $X,$
defined by C. Cibils and M. Rosso in~\cite{CR}, and view all Hopf
bimodules as left $X$-modules. The above cohomologies will then be isomorphic to the space of extensions $Ext_X^* (M,N).$

We will further give expressions of the Yoneda product of extensions in $H_{A4}^* (M,N )$ (which generalizes $H_{A4}^*(M,A))$  and $H^*_b(A,A),$ which do not involve $X.$

There appears to be a natural definition for the cyclic cohomology of a finite dimensional Hopf algebra, that is the cyclic cohomology of the associative algebra $X.$ On the other hand, in~\cite{CM}, A. Connes and H. Moscovici have defined a cyclic module associated to a Hopf algebra endowed with a modular pair in involution.  These theories differ; moreover, although the cyclic cohomology of $X$ takes equally into account all the structures of $A,$ the cyclic cohomology of $A$ defined by A. Connes and H. Moscovici  gives a greater weight to the coalgebra structure of $A.$ However, it would be interesting to see whether there is a connection between them.

This paper is organized as follows: we will first describe the algebra $X.$ Then we will describe M. Gerstenhaber and S.D. Schack's cohomologies $H_{GS}^*(M,N)$ and $H^*_b(A,A)$, explain why the first extends the second, and prove that $H_{GS}^*(M,N)$ is isomorphic to $Ext_X^* (M,N).$
Next, we shall do the same with a generalization of the cohomology defined by C. Ospel. Finally, we shall describe a cup-product on this last cohomology and on $H^*_b(A,A),$ and prove that it corresponds to the Yoneda product of extensions.

We shall need some notation: all tensor products will be taken over
the base field $k.$ The structure maps of the Hopf algebra $A$ are
denoted in the standard way $\mu, \ \Delta, \  \eta,  \ \varepsilon$
and $S.$ We will use Sweedler's notations for the comultiplications: 
\begin{eqnarray*}
\Delta (x) & = & \sum x^{(1)} \otimes x^{(2)} \mathrm{\ if \ } x \in A, \\
\mu ^* (l) & = & \sum l_{(1)} \otimes l_{(2)} \mathrm{\ if \ } l \in A^*,
\end{eqnarray*} and the comodule structures:
$$\begin{array}{ccccccc}
M  & \longrightarrow & A \otimes M & \mbox{ and } & M  & \longrightarrow & M \otimes A \\
m & \mapsto & \sum m_{(-1)} \otimes m_{(0)} & & m & \mapsto & \sum m_{(0)} \otimes m_{(1)}.
\end{array}$$ 
Let us introduce, for any integer $u \geq -1,$ the map:
\begin{eqnarray*}
\Delta ^{(u)} : A & \longrightarrow & A^{\otimes u+1} \\
a & \mapsto & \begin{cases}
\varepsilon (a) & \text{if $u=-1$}, \\
a  & \text{if $u=0$},\\
(\Delta \otimes 1^{\otimes u-1} ) \ldots (\Delta \otimes 1) \Delta (a) & \text{otherwise.} \\

\end{cases}
\end{eqnarray*}

The results in this paper will form part of my PhD thesis. I am happy to thank my advisor Claude Cibils and Jean-Michel Oudom for their help and comments.


\section{The algebra $X.$ }

Let $A$ be a finite dimensional Hopf algebra over a field $k.$ In~\cite{CR}, C. Cibils and M. Rosso have introduced an associative $k-$algebra $X$, in order to identify Hopf bimodules over $A$ with left modules over $X$ (just as bimodules over an associative algebra are identified with left modules over its enveloping algebra). This they do by viewing the right $A-$module structure of a Hopf bimodule as a left $A^{op}-$module stucture, its left coaction as a left action of $A^{*op},$ and its right coaction as a left action of $A^*.$ To preserve the compatibilities between the structures, they endow $A^*$ and $A^{*op}$ with $A-$bimodule structures involving the antipode. Thus they define the algebra $X$ as being the vector space $(A^{*op} \otimes A^* ) \underline{ \otimes} (A \otimes A^{op})$ in which the first two and last two tensorands keep natural multiplication, and 
\begin{multline}
[(1 \otimes 1)  \underline{ \otimes} (a \otimes b)] \,  [ (l \otimes k) \underline{  \otimes} (1 \otimes 1)]\\ 
= \sum l_{(1)} (Sa^{(1)}) k_{(1)} (S^{-1}a^{(3)}) 
 l_{(3)} (S^{-1}b^{(1)}) k_{(3)} (Sb^{(3)}) \  [ (l_{(2)} \otimes k_{(2)} )
\underline{ \otimes} ( a^{(2)} \otimes b^{(2)}) ],
\notag
\end{multline} and they prove the following:

\begin{theo} \textnormal{(\cite{CR} theorem 3.1)} 
There is a vector space-preserving equivalence of categories between the category of  left modules over $X$ and the category of Hopf bimodules over $A.$ 

\end{theo} 

\paragraph{Remarks} (i) If $M$ is a Hopf bimodule, its left $X-$module structure is given by:
$$ [(l \otimes k)\underline{ \otimes} (a \otimes b) ].m = \sum l( a^{(1)} m_{(-1)}  b^{(1)} ) k( a^{(3)} m_{(1)}  b^{(3)} ) ( a^{(2)} m_{(0)}  b^{(2)}). $$

(ii) C. Cibils and M. Rosso have also proved that the algebra $X$ is isomorphic to the direct tensor product $\mathcal{H} (A) \otimes \mathcal{D} (A) ^{op},$ where $\mathcal{H} (A)$ is the Heisenberg double of $A,$ and $\mathcal{D} (A)$ is the Drinfel'd double of $A.$

\

Henceforth, we will identify Hopf bimodules over $A$ with left $X-$modules, and, whenever $M$ and $N$ are Hopf bimodules, we shall consider the algebra  $Ext_X^* (M,N).$


\section{Gerstenhaber and Schack's cohomology for Hopf bimodules}

This has been defined by M. Gerstenhaber
 and S.D. Schack in~\cite{GS2}; we shall also use various notions introduced by S. Shnider and S. Sternberg in~\cite{Sh-St}.

\subsection{Bar resolution of a Hopf bimodule}

Let us first consider a  bimodule $M$ over $A.$ Recall that the bar complex $B _{\bullet} (M)$ of $M$ is the direct sum $\bigoplus_{q \geq -1} B_q (M),$ in which the $B_q (M)=A^{\otimes q+1} \otimes M \otimes A^{\otimes q+1} $ are bimodules via the standard actions (\emph{i.e.} multiplication of the leftmost or rightmost tensorand), the differential $\partial _q \, : \, B_q (M)  \rightarrow   B_{q-1} (M)$ sends
$ a_0 \otimes \ldots \otimes a_q \otimes m \otimes b_q \otimes \ldots \otimes b_0 $ to 
\begin{multline}  \sum _{i=0} ^{q-1} (-1)^i a_0 \otimes \ldots \otimes a_i a_{i+1} \otimes \ldots \otimes m \otimes \ldots \otimes b_{i+1} b_i \otimes \ldots \otimes b_0\\
+ (-1)^q a_0 \otimes \ldots \otimes a_{q-1} \otimes a_q m b_q \otimes b_{q-1} \otimes \ldots \otimes b_0 
\notag
\end{multline}
for $q \geq 0,$ and $\partial _{-1} =0.$ This defines a complex of $A-$bimodules. To see that it is a resolution of $M,$ let us introduce the homotopy

$$\begin{array}{lccl}
h_q \;  : & B_q (M) & \rightarrow  & B_{q+1} (M) \\
 & x & \mapsto & 1 \otimes x \otimes 1.
\end{array}$$
It satisfies the equation $\partial _{q+1} h_q + h_{q-1} \partial _q = id,$ for all $q \geq -1.$

Assume now that $M$ is a Hopf bimodule over $A.$ Let us endow the $B_q (M)$ with codiagonal coactions, that is \begin{eqnarray*}
B_q (M) & \rightarrow  & A \otimes B_{q} (M) \\
\underline{a} \otimes m \otimes \underline{b} & \mapsto &  \sum a^{(1)}_0 \ldots a^{(1)}_q m_{(-1)} b^{(1)}_q \ldots b^{(1)}_0 \\
& & \ \otimes a^{(2)}_0 \otimes \ldots \otimes  a^{(2)}_q \otimes  m_{(0)} \otimes   b^{(2)}_q \otimes \ldots \otimes  b^{(2)}_0, 
\end{eqnarray*}
and similarly on the right. This makes $B _{\bullet} (M)$ into a complex of Hopf bimodules. The homotopy becomes a morphism of  $A-$bicomodules, and $B_{\bullet} (M)$ is therefore a Hopf bimodule resolution of $M$ which splits as a sequence of $A-$bicomodules.

\subsection{Cobar resolution of a Hopf bimodule}

Dually, let $N$ be a bicomodule over $A.$ The cobar complex $C^
{\bullet} (N) $ is equal to  the direct sum $\bigoplus_{p \geq -1} C^p (N),$ where the vector spaces $C^p (N)=A^{\otimes p+1} \otimes N \otimes A^{\otimes p+1} $ are bicomodules via the standard coactions (that is comutiplication on the first or last tensorand), the differential $\partial ^p \, : \, C^p (N)  \rightarrow   C^{p+1} (N)$ sends\\ $ a_0 \otimes \ldots \otimes a_p \otimes n \otimes b_p \otimes \ldots \otimes b_0 $ to 
\begin{multline}  \sum _{i=0} ^{p} (-1)^i a_0 \otimes \ldots \otimes \Delta a_i  \otimes \ldots \otimes n \otimes \ldots \otimes \Delta b_i \otimes \ldots \otimes b_0\\
+ (-1)^{p+1} \sum a_0 \otimes \ldots \otimes a_{p} \otimes  n_{(-1)} \otimes n_{(0)} \otimes n_{(1)} \otimes b_p \otimes \ldots \otimes b_0, 
\notag
\end{multline} 
and $\partial^{-1} =0.$

There is a homotopy $$\begin{array}{rcl}
h^p \;  :  C^p (N) & \rightarrow  & C^{p-1} (N) \\
\underline {a} \otimes  n \otimes \underline{b} & \mapsto & \varepsilon (a_0) \varepsilon (b_0) a_1 \otimes \ldots \otimes a_p \otimes n \otimes b_p \otimes \ldots \otimes b_1 
\end{array} $$ which enables us to see that $C^{\bullet} (N)$ is a bicomodule resolution of $N.$

When $N$ is a Hopf bimodule, the $C^p (N)$ are endowed with diagonal actions, that is \begin{eqnarray*}
A \otimes C^p (N) & \rightarrow  & C^p (N) \\
c \otimes \underline {a} \otimes  n \otimes \underline{b}& \mapsto & \sum c^{(1)} a_0 \otimes \ldots \otimes c^{(p+1)} a_p \otimes c^{(p+2)} n \otimes  \ldots \otimes c^{(2p+3)} b_0 \\
& &\  = \Delta ^{(2p+2)} (c) \;  (a_0 \otimes \ldots \otimes a_p \otimes n \otimes b_p \otimes \ldots \otimes b_0), 
\end{eqnarray*} and similarly on the right, so that $C^ {\bullet} (N)$ becomes a Hopf bimodule resolution of $N$ which splits as a sequence of $A-$bimodules.

\subsection{Gerstenhaber and Schack cohomology of Hopf bimodules}

Let $M$ and $N$ be two Hopf bimodules over $A.$ As in~\cite{GS2}, set $$C^{pq} (M,N)= Hom_{A4} (B_q (M),C^p (N))$$
 (that is the space of Hopf bimodule morphisms from $B_q (M)$ to $C^p (N)$), 
 
$$\begin{array}{rcl}
\delta _h \; : C^{pq} (M,N) & \rightarrow & C^{p,q+1} (M,N) 
\\
\alpha & \mapsto & \alpha \circ \partial _{q+1}
\end{array}$$
and $$\begin{array}{rcl}
\delta _c \; : C^{pq} (M,N) & \rightarrow & C^{p+1,q} (M,N) 
\\
\alpha & \mapsto & (-1)^q  \partial ^{p} \circ \alpha .
\end{array}$$

This defines a double complex, and the homology of its total complex is Gerstenhaber and Schack's cohomology of $M$ into $N;$ let $H^* _{GS} (M,N)$ denote this.

\paragraph{Remark}
This double complex is isomorphic to the double complex whose modules are $$Hom_k (A^{ \otimes q} \otimes M \otimes A^{ \otimes q}, A^{ \otimes p} \otimes N \otimes A^{ \otimes p} ),$$ with differentials $\delta '_h$ and $\delta '_c$ defined by: \begin{align} 
\delta '_h ( \alpha )(\mathbf{a_{0,q}} \otimes m \otimes
\mathbf{b_{q,0}}) = & \ \Delta ^{(2q)} (a_0) \,  \alpha ( \mathbf{a_{1,q}} \otimes  m \otimes \mathbf{b_{q,1}}) \, \Delta ^{(2q)}(b_0) \notag \\
\  & + \sum _{i=1} ^{q} (-1)^i \, \alpha (a_0 \otimes \ldots \otimes a_{i-1} a_{i} \otimes \ldots \otimes b_{i} b_{i-1} \otimes \ldots )  \notag \\
\ & + (-1)^{q+1} \, \alpha (a_0 \otimes \ldots \otimes a_q m b_q \otimes \ldots \otimes b_0) \notag
\end{align}
and \begin{align} 
(-1)^q\delta '_c ( \alpha ) = & (1 \otimes \alpha \ot 1) \circ (\delta
_L \ot 1)  \circ \delta _R+ \sum _{i=1} ^{p}(-1)^i \tilde \Delta _i \circ \alpha \notag \\
\ & + (-1)^{p+1} [1 ^{\otimes p} \otimes ((1  \otimes \delta _R ^{(N)}) \circ
\delta _L ^
{(N)} ) \otimes  1 ^{\otimes p } ] \circ \alpha, \notag
\end{align} where $\tilde \Delta _i = \Delta _{2p-i+1} \circ \Delta
_i,$ the notation   $\mathbf{a_{0,q}}$ stands for $a_0 \ot \ldots \ot
a_q,$ and $ \mathbf{b_{q,0}}=b_q \ot \ldots \ot b_0. $ The various
$\delta$ are the coactions on the appropriate Hopf bimodules.

\subsection{Gerstenhaber and Schack cohomology of a Hopf algebra}\label{sect2.4}

When $M$ and $N$ are both equal to the Hopf bimodule $A,$
M. Gerstenhaber and S.D. Schack have given, in a similar way, a definition for the cohomology of $A$ (see~\cite{GS1}):

Consider the bar resolution $Bar_{\bullet} (A)$ of $A \ (Bar_q (A) =
A^{\otimes q+2}),$ and the cobar resolution $Cob ^{\bullet} (A)$ of $A
\ (Cob^p (A) = A^{\otimes p+2}):$ they enjoy the same properties as
the resolutions $B_{\bullet} (A)$ and $C^{\bullet} (A),$ and the
double complex associated to these resolutions is a double complex
whose modules are isomorphic to $C ^{\prime pq} =Hom_k (A^{\otimes q}
, A^{\otimes p}).$ The homology of its total complex is denoted by $H_b^*(A,A).$ We shall see further on that $H_{GS}^*$ generalizes $H_b^*,$ that is  $H_{GS}^* (A,A)\cong H_b^*(A,A).$

\paragraph{Example} When $A$ is a  group bialgebra $kG,$ it is cosemisimple, so that the rows of the double complex are acyclic, and the homology of the double complex is equal to the homology of the remaining column, that is the cohomology  $H^*(G,k)=HH^*(kG,k)$ (see~\cite{GS1}. Further, various cases of group bialgebras have been studied in~\cite{PW}).

\paragraph{Remark} There is a suitable sub-double complex $\hat{C} ^{\bullet \bullet}$ of $ C ^{\prime \bullet \bullet}$ whose homology is adapted to the study of deformations of the Hopf algebra $A$ as a Hopf algebra. It is obtained from  $ C ^{\prime \bullet \bullet}$ by replacing the first row and the first column by $0 \rightarrow 0 \rightarrow \ldots 
\rightarrow 0 \rightarrow \ldots$ (see~\cite{GS1} p79). It has been used by A. Giaquinto in his thesis (\cite{G}) to study the preferred deformations of the classical matrix bialgebra and plane (that is deformations in which the comultiplication remains unchanged); he interprets the quantum matrix bialgebra $\mathbf{\mathrm M}_q(2)$ and the quantum plane $\mathbf{\mathrm k}_q^2$ as such.


\subsection{Link with extensions}

In this paragraph, we shall  prove that Gerstenhaber and Schack's cohomology $H^* _{GS} (M,N),$ and $Ext ^* _X (M,N),$ are the same. We will first need some definitions and a few preliminary results.

\paragraph{Definition}
A Hopf bimodule M is said to be a \textit{relative $A$-4 projective} when the functor $Hom_{A4} (M,-)$ associates exact sequences of vector spaces to exact sequences of Hopf bimodules which split as sequences of $A-$bicomodules. A resolution of a Hopf bimodule is called a \textit{ relative $A$-4 projective resolution} if all its terms are relative $A$-4  projectives and if the sequence splits as  a sequence of $A-$bicomodules. One defines dually relative $A$-4 injectives  and relative $A$-4 injective resolutions (see~\cite{Sh-St} chapter 10).

\paragraph{Example} A projective Hopf bimodule is a relative $A$-4
projective, an injective Hopf bimodule is  a relative $A$-4 injective.

\paragraph{Example} The resolutions $B_{\bullet} (M)$ and   $Bar_{\bullet} (A)$ are relative $A$-4 projective resolutions; the resolutions $C^{\bullet} (N)$ and $Cob^{\bullet} (A)$ are relative $A$-4 injective resolutions. 

Indeed, let us for instance look at $B_{\bullet} (M).$ In the first place, for any Hopf bimodule $Y,$ there is an isomorphism between $Hom_{A4}(B_{q} (M), Y)$ and the space of bicomodule morphisms $Hom^{AA}(A^{\otimes q} \otimes M \otimes A^{\otimes q}, Y)$, and the functor $Hom^{AA}(A^{\otimes q} \otimes M \otimes A^{\otimes q},-)$ clearly takes exact sequences of Hopf bimodules which split as sequences of $A-$bicomodules to exact sequences of vector spaces, so that each $B_{q} (M)$ is a relative $A$-4 projective. Furthermore, we have said earlier that the Bar resolution splits as a sequence of $A-$bicomodules, so that $B_{\bullet} (M)$  is a  relative $A$-4 projective resolution of $M.$ Similar arguments can be applied to the other resolutions.

\   

These relative resolutions have properties which are similar to those of projective resolutions:
 
\begin{prop}\label{projequiv} Let $Y$ and $Z$ be two relative $A$-4 projective  resolutions of a Hopf bimodule M. Then there exists a map of complexes of Hopf bimodules $\phi$ such that the triangle $$\xymatrix{ Y \ar[rr]^{\phi} \ar[rd] & & Z \ar[ld]  \\
& M & }$$ commutes, and $\phi$ is unique up to homotopy of Hopf
bimodule complexes. A similar statement is true for relative $A$-4 injective resolutions.
\end{prop}
\textit{Proof}: Follow the proof of the fact that two projective resolutions are homotopic (see~\cite{Sh-St} p242-245), bearing in mind that all the exact sequences involved split as sequences of $A-$bicomodules. $\square$

\begin{cor} Let $Y$ and $T$ be two relative  $A$-4 projective resolutions of a Hopf bimodule $M,$ and let $Z$ and $U$ be two relative $A$-4 injective resolutions of a Hopf bimodule $N.$ Then the double complexes $Hom_{A4} (Y,Z)$ and $Hom_{A4} (T,U)$ endowed with the natural differentials are homotopy equivalent.
\end{cor}
\textit{Proof}: See~\cite{Sh-St} p246. $\square$

\paragraph{Consequences} (i) Recall that $B_{\bullet} (M)$ is a
relative  $A$-4 projective resolution of $M,$ therefore in the
computation of $H_{GS} ^* (M,N),$ the resolution $B_{\bullet} (M)$ can
be replaced by any relative $A$-4 projective resolution of $M$. The
resolution $C^{\bullet}(N)$ can also be replaced by any relative $A$-4
injective resolution of $N.$

(ii) In particular,  $H_{GS} ^* (M,N)$ generalizes $H_b ^* (A,A).$ 

\ 

We can now prove:

\begin{theo}\label{gsisom} Let $M$ and $N$ be  Hopf bimodules over $A.$ There is an isomorphism $H_{GS} ^* (M,N) \cong Ext_X ^* (M,N).$
\end{theo} 
\textit{Proof}: Let us consider the functor $Ext_X^*(-,N),$ with $N$
fixed. It is characterized by the following properties (\emph{cf.}
\cite{McL}):
\begin{enumerate}
\item $Ext_X^0(M,N) \cong Hom_X(M,N)=Hom_{A4}(M,N),$
\item $Ext_X^n(P,N)=0$ for every $n\geq 1$ and every projective Hopf
bimodule $P,$
\item $Ext_X^*(-,N)$ is a cohomological $\delta$-functor (\emph{see}
\cite{W} p30). \end{enumerate}

It is therefore enough to check that the functor $H_{GS}^*(-,N)$
satisfies these properties. So the theorem follows from the following
three lemmas:

\begin{lemma}\label{lemma1} For every Hopf bimodule $M,$ there is an
isomorphism $$H^0_{GS}(M,N) \cong Hom_{A4}(M,N).$$ \end{lemma}

Let $\alpha$ be an element in  $Hom _{A4}(A \ot M\ot A,A \ot
N \ot A).$ It is a cocycle if and only if $\delta_h(\alpha)=0$ and
$\delta_c(\alpha)=0,$ that is $\alpha \circ
\partial_1=0=\partial^1 \circ \alpha ,$ or $\alpha$ is zero on ${\rm
Im}\, \partial _1 =\ker \partial _0$ and $\alpha $ takes values in $\ker
\partial ^1 = {\rm Im} \, \partial  ^0,$ with $\partial _0$ onto
and $\partial ^0$ one-to-one. It follows that if $\alpha$ is a
cocycle, it yields a unique morphism $\overline{\alpha} \in
Hom _{A4}(M,N)$ such that $\partial^0 \circ \overline{\alpha} \circ \partial
_0=\alpha.$ On the other hand, if $\beta$ is an element of
$Hom_{A4}(M,N),$ then $\partial^0 \circ \beta \circ \partial
_0$ is a cocycle.

\begin{lemma}\label{lemma2} For every projective Hopf bimodule $P$ and
every integer $n \geq 1,$ the $k-$vector space $H_{GS}^n(P,N)$
vanishes.\end{lemma}

Let $P$ be a projective Hopf bimodule. Then $\cdots 0 \rightarrow 0
\rightarrow \cdots \rightarrow 0 \rightarrow P \rightarrow P
\rightarrow 0$ is a relative $A$-4 projective resolution of $P.$ It
is therefore homotopy equivalent to the bar resolution $B_{\bullet}(M)$
(proposition \ref{projequiv}). Gerstenhaber and Schack's cohomology is
therefore the homology of the double complex in which all the terms
are zero except those on the first line, which is equal to
$Hom_{A4}(P,C^{\bullet}(N)).$ This line is acyclic, since $P$ is
projective and $C^{\bullet}(N)$ is exact. Its cohomology is therefore
zero in positive degree. 

\begin{lemma}\label{lemma3} The functor $H_{GS}^*(-,N)$ takes short
exact sequences of Hopf bimodules $0 \rightarrow M' \rightarrow M
\rightarrow M'' \rightarrow 0$ to long exact sequences of vector
spaces $$\cdots \rightarrow H_{GS}^n(M'',N)\rightarrow
H_{GS}^n(M,N)\rightarrow H_{GS}^n(M',N)\stackrel{\delta}{\rightarrow}
H_{GS}^{n+1}(M'',N)\rightarrow \cdots$$ \end{lemma}

Let $0 \rightarrow M' \rightarrow M
\rightarrow M'' \rightarrow 0$ be an exact sequence of Hopf
bimodules. Then, for every $q\geq 0,$ the sequence $$(E_q): 0
\rightarrow A^{\ot q}\ot M' \ot A ^{\ot q} \rightarrow A^{\ot q}\ot M
\ot A ^{\ot q} \rightarrow A^{\ot q}\ot M'' \ot A ^{\ot q}
\rightarrow0$$ is an exact sequence of vector spaces. The sequences 
\begin{eqnarray*}0 \rightarrow && Hom _k (A^{\otimes q} \ot M'' \otimes
A^{\otimes q} , A^{\otimes p}  \otimes N \ot  A^{\otimes p}  )\\ && \rightarrow  Hom _k (A^{\otimes q} \ot M \otimes A^{\otimes q} , A^{\otimes p}  \otimes N \ot  A^{\otimes p} ) \\
 &&  \rightarrow   Hom _k (A^{\otimes q} \ot M' \otimes A^{\otimes q} , A^{\otimes p}  \otimes N \ot  A^{\otimes p} ) \rightarrow 0
\end{eqnarray*} are therefore exact for all integers $p, q \geq 0.$
These exact sequences commute with the differentials, and thus define
an exact sequence of double complexes. This yields a long exact sequence $$\cdots \rightarrow H_{GS}^n(M'',N)\rightarrow
H_{GS}^n(M,N)\rightarrow H_{GS}^n(M',N)\stackrel{\delta}{\rightarrow}
H_{GS}^{n+1}(M'',N)\rightarrow \cdots \hspace{1cm} \square $$

\paragraph{Remark} To compute $ H_{GS}^*(M,N) \cong Ext_X ^* (M,N),$ one can also take a projective resolution $P _{\bullet}$ of $M$ (or an injective resolution $I^{\bullet}$ of $N$) and take the homology of the complex $Hom_{A4} (P _{\bullet} ,N) $ (or of the complex $Hom_{A4} (M,I^{\bullet})$).


\section{ Hopf bimodule cohomology}

When $A$ is a  Hopf algebra and $M$ is a Hopf bimodule over $A,$ C. Ospel has defined in his thesis~\cite{O} a cohomology of $M$ with coefficients in $A.$ This can be extended in the following way:

Let $M$ and $N$ be Hopf bimodules. Consider the vector spaces $M \otimes A^{\otimes {q+1}}$ and $A^{\otimes {p+1}} \otimes N$ for all non negative integers $p$ and $q.$ Each of these vector spaces is endowed with a Hopf bimodule structure:

The space $M \otimes A^{\otimes {q+1}}$ is equipped with the standard actions and codiagonal coactions, that is $$\begin{array}{cccl}
\mu _L \, : & A \otimes  M \otimes A^{\otimes {q+1}} & \rightarrow & M \otimes A^{\otimes {q+1}} \\
\  & b \otimes m \otimes \underline{a} & \mapsto & bm \otimes \underline{a} \\
\delta _L \, : & M \otimes A^{\otimes {q+1}} & \rightarrow & A \otimes  M \otimes A^{\otimes {q+1}} \\
\  & m \otimes \underline{a}  & \mapsto & \sum m_{(-1)} a^{(1)} _1 \ldots a^{(1)} _{q+1} \otimes m_{(0)} \otimes  a^{(2)} _1 \otimes \ldots \otimes a^{(2)} _{q+1} 
\end{array}$$ and similarly on the right. Dually,  $A^{\otimes {p+1}} \otimes N$ is endowed with the diagonal actions and standard coactions : $$\begin{array}{cccl}
\mu _L \, : & A \otimes A^{\otimes {p+1}} \otimes N & \rightarrow &  A^{\otimes {p+1}} \otimes N \\
\ & b \otimes \underline{a} \otimes n & \mapsto &  \sum b^{(1)} a_1 \otimes \ldots \otimes b^{({p+1})} a_{p+1} \otimes b^{(p+2)} n \\
\  & \  & \ & \  = \Delta ^{({p+1})} (b) \, ( \underline{a} \otimes n) \\ 
\\
\delta _L \, : &  A^{\otimes {p+1}} \otimes N & \rightarrow &  A \otimes A^{\otimes {p+1}} \otimes N \\
\  & \underline{a} \otimes n & \mapsto &  \sum a^{(1)} _1 \otimes  a^{(1)} _2 \otimes  a_2 \otimes \ldots \otimes  a_{p+1} \otimes n,
\end{array}$$ with similar action and coaction on the right.

For all non negative integers $p$ and $q,$ set $$K^{pq} (M,N) = Hom _{A4} (M \otimes A^{\otimes q+1},  A^{\otimes p+1} \otimes N),$$ that is the space of Hopf bimodule morphisms from $M \otimes A^{\otimes q+1}$ to $ A^{\otimes p+1} \otimes N,$ and consider the maps
$$\begin{array}{ccl}
d_h \, : \, K^{pq} (M,N) & \rightarrow & K^{p,q+1} (M,N) \\
\alpha & \mapsto & \alpha \circ \lambda \\
d_c \, : \, K^{pq} (M,N) & \rightarrow & K^{p+1,q} (M,N) \\
\alpha & \mapsto & (-1)^q \rho \circ \alpha,
\end{array}$$ where $$ \begin{array}{ccl}
\lambda \, : \,  M \otimes A^{\otimes q+2} & \rightarrow &  M \otimes A^{\otimes q+1} \\
\ m \otimes \underline{a}  & \mapsto & m \, a_0 \otimes \ldots \otimes a_{q+1} \\
\ & \ & \  + \sum _{i=0} ^{q} (-1)^{i+1} m \otimes \ldots \otimes a_i a_{i+1} \otimes \ldots \otimes a_{q+1},\  \mathrm{and} \\
 & &  \\
\rho \, : \,  A^{\otimes p+1} \otimes N & \rightarrow &  A^{\otimes p+2} \otimes N \\
\ \ \ x  & \mapsto & (\sum _{i=0} ^{p} (-1)^i \Delta_i + (-1)^{p+1} (1 ^{\otimes p} \otimes \delta _L^N))(x).
\end{array}$$

This defines a double complex. Let $H_{A4} ^* (M,N)$ denote the homology of its total complex, which we shall simply call \textit{Hopf bimodule cohomology}. 

Note that $(M\ot A^{\ot \bullet +1},\lambda)$ is the bar resolution of
$M$ as a \emph{right} module, and $(A^{\ot \bullet +1}\ot N,\rho)$ is
the cobar resolution of $N$ as a \emph{left} comodule.

\paragraph{Remarks}
(i) This double complex is isomorphic to the double complex whose modules are the spaces of left $A-$module and right $A-$comodule maps $$Hom _{A-}^{-A} ( M \otimes A^{\otimes q}, A^{\otimes p} \otimes N),$$ the vertical differential translating as $\alpha \mapsto d(\alpha)$ with  \begin{eqnarray*}
d(\alpha) (m \otimes a_1 \otimes \ldots \otimes a_{q+1}) & = & \alpha ( m \, a_1 \otimes \ldots \otimes a_{q+1}) \\
\ & & +\sum _{i=1} ^{q} (-1)^{i} \alpha (m \otimes \ldots \otimes a_i a_{i+1} \otimes \ldots \otimes a_{q+1}) \\
\ & & + (-1)^{q+1} \alpha ( m \otimes a_1  \otimes \ldots \otimes a_q ).a_{q+1}, 
\end{eqnarray*} and the horizontal differential as $$\alpha \mapsto (1 \otimes \alpha ) \circ \delta _L + \sum _{i=1}^{p} (-1)^i \Delta _i \circ \alpha + (-1)^{p+1} (1 ^{\otimes p } \otimes \delta _L ^N ) \circ \alpha.$$

Furthermore, when $N$ is the Hopf bimodule $A, \  K^{\bullet \bullet } (M,A)$ is isomorphic to the double complex defined by C. Ospel in his thesis~\cite{O}.

ii)  When $M$ and $N$ are both equal to  $A,\   K^{\bullet \bullet } (A,A)$ is isomorphic to the double complex which defines Gerstenhaber and Schack's bialgebra cohomology $H_b^* (A,A)$ (see~\cite{GS1} and subsection~\ref{sect2.4}): $H_{A4}^* (M,N)$ generalizes $H_b^*(A,A).$

(iii) When $N$ is equal to $A, $ the homology of the $p^{th}$ column is the Hochschild cohomology of $A$ with coefficients in $Hom_{A-} (M,A ^{\otimes p} ),$ the space of left $A-$module morphisms from $M$ to $A ^{\otimes p},$ this last space being endowed with an $A-$bimodule structure as follows: $A ^{\otimes p}$ is equipped with the diagonal $A-$bimodule structure, and if $f \in Hom_{A-} (M,A ^{\otimes p}), \ a \in A$ and $v \in M,$ $$(a.f)(v)=f(v.a) \mbox{ and } (f.a)(v) = f(v)a.$$ 
The homology of the $q^{th}$ row is the Cartier cohomology $H_c^* (M \otimes  A ^{\otimes q} , A)$ (see~\cite{C}).

We shall now study the relationship between this cohomology and the space $Ext _X^* (M,N).$

First, let us state a slight extension of a  theorem proved by
C. Ospel. We shall also provide a short sketch of the proof, as it will be useful in the proof of theorem \ref{Yonedath}.

\begin{theo}\label{Ospel} \textup{(see \cite{O} theorem 3.1)} Let $M$ and $N$ be  Hopf bimodules. There is an isomorphism $$H_{A4}^1 (M,N) \cong Ext_X^1 (M,N).$$
\end{theo}
\textit{Proof}: Take a 1-cochain $f=(f_0,f_1)$  in $Hom _{A-}^{-A} (M \otimes A, N) \oplus Hom _{A-}^{-A} (M ,A \otimes  N).$ The corresponding 1-extension $F$ is equal to $$0 \rightarrow N \rightarrow N \oplus M \rightarrow M \rightarrow 0$$ as a sequence of vector spaces, and the $A-$bimodule structure on  $ N \oplus M$ is given by:
$$\begin{pmatrix} \mu _R ^N & f_0 \\ 0 & \mu _R ^M \end{pmatrix} \mathrm{ \ for\  the \  right \  A-module \  structure,} $$
$$\begin{pmatrix} \delta _L ^N & f_1 \\ 0 & \delta _L ^M \end{pmatrix} \mathrm{ \ for \  the \  left \  A-comodule \  structure,}$$ $$\mathrm{and \ } \begin{pmatrix} \mu _L ^N & 0 \\ 0 & \mu _L ^M \end{pmatrix} \ \mathrm{ and} \ \begin{pmatrix} \delta _R ^N & 0 \\ 0 & \delta _R ^M \end{pmatrix} \mathrm{ \  for \  the \  remaining \  structures.}$$

The sequence $F$ is a sequence of Hopf bimodules if and only if  $f$ is a 1-cocycle. This describes the isomorphism. $\square$

\ 

More generally, the following is true:

\begin{theo}\label{Ospelth} For any  Hopf bimodules $M$ and $N,$ there is an isomorphism $$H_{A4}^*(M,N) \cong Ext _X^* (M,N).$$
\end{theo}
\textit{Proof}: As in theorem \ref{gsisom}, we shall prove that the
functor $H_{A4}^*(-,N)$ satisfies the properties of the functor $Ext _X^* (-,N): $

\begin{lemma} For every Hopf bimodule $M,$ there is an
isomorphism  $$H_{A4}^0(M,N) \cong Hom_{A4} (M,N).$$ 
\end{lemma}

This is proved as in lemma \ref{lemma1}.

\begin{lemma} For every projective Hopf bimodule  $P$ and every
integer $n \geq 1,$ the $k-$vector space $H_{A4}^n(P,N)$ vanishes.
\end{lemma}

Let $P$ be a projective Hopf bimodule. Let us consider its bar
resolution $$B_{\bullet}^A(P): \cdots \rightarrow P \ot A^{\ot q+2}
\stackrel{\lambda_q}{\longrightarrow}P \ot A^{\ot q+1}
\stackrel{\lambda_{q-1}}{\longrightarrow}\cdots \rightarrow P \ot
A^{\ot 2} \stackrel{\lambda_0}{\longrightarrow}P \ot A
\stackrel{\lambda_{-1}}{\longrightarrow}P \rightarrow 0,$$ with
$\lambda _{-1}(u \ot a)=u.a.$

Since $P$ is projective and $\lambda _{-1}$ is onto, $\lambda _{-1}$
has a section: there exists a morphism of Hopf bimodules $s: P
\rightarrow P \ot A$ satisfying $\lambda _{-1}(s(u))=u$ for every $u
\in P.$

Set  \begin{eqnarray*} h_q\ :P\ot A^{\ot q+1} & \longrightarrow &
P\ot A^{\ot q+2}\\ u\ot a_0 \ot \ldots \ot a_q &\mapsto & s(u)\ot
a_0 \ot \ldots \ot a_q.\end{eqnarray*}
It is a Hopf bimodule morphism (recall that the Hopf bimodules are
endowed with standard actions and codiagonal coactions). Set
$s(u)=\sum_{i=1}^{t}v_i \ot b_i,$ and let us compute $\lambda_q h_q+h_{q-1}\lambda_{q-1}:$
\begin{eqnarray*}
\lambda _q h_q (u\ot a_0 \ot \ldots \ot a_q)&=& \lambda _q (s(u)\ot
a_0 \ot \ldots \ot a_q)\\&=& \lambda _q ( \sum_{i=1}^{t}v_i \ot
b_i \ot
a_0 \ot \ldots \ot a_q)\\
&=& \sum_{i=1}^{t} [ v_i. b_i \ot a_0 \ot \ldots \ot a_q \\
&& \hspace{.8cm}- v_i \ot b_i a_0 \ot a_1 \ot \ldots \ot a_q \\
&& \hspace{.8cm}+\sum_{j=0}^{q-1}(-1)^j\,  v_i \ot b_i \ot a_0 \ot \ldots
\ot a_j a_{j+1}\ot \ldots \ot a_q
 ]\\
&=& \lambda _{-1}(s(u))\ot a_0 \ot \ldots \ot a_q\\
&& \hspace{.8cm}-(s(u).a_0) \ot a_1 \ldots \ot a_q \\
&& \hspace{.8cm}+  \sum_{j=0}^{q-1}(-1)^j \, s(u)\ot a_0 \ot \ldots
\ot a_j a_{j+1}\ot \ldots \ot a_q\\
&=& u \ot a_0 \ot \ldots \ot a_q \\
&& \hspace{.8cm}- s(u.a_0) \ot a_1 \ot \ldots \ot a_q\\
&& \hspace{.8cm}+ \sum_{j=0}^{q-1}(-1)^j\, s(u) \ot a_0 \ot \ldots
\ot a_j a_{j+1}\ot \ldots \ot a_q
 \end{eqnarray*}   

and \begin{eqnarray*} 
h_{q-1}\lambda _{q-1}(u\ot a_0 \ot \ldots \ot a_q)&=&
h_{q-1}(u.a_0 \ot a_1 \ot \ldots \ot a_q)\\
&& +  \sum_{j=0}^{q-1}(-1)^{j+1}\, h_{q-1}(u \ot a_0 \ot \ldots
\ot a_j a_{j+1}\ot \ldots \ot a_q)\\
&=& s(u.a_0) \ot a_1 \ot \ldots \ot a_q\\
&& -   \sum_{j=0}^{q-1}(-1)^j\, s(u) \ot a_0 \ot \ldots
\ot a_j a_{j+1}\ot \ldots \ot a_q.
\end{eqnarray*}

Adding these expressions yields: $$(\lambda _q
h_q + h_{q-1}\lambda _{q-1})(u\ot a_0 \ot \ldots \ot a_q)=u\ot
a_0 \ot \ldots \ot a_q.$$ Therefore $h_{\bullet} $ is a Hopf bimodule
homotopy from $id$ to $0.$

Now fix $p \in \mathbb{N},$ and consider the complex $Hom_{A4}(P \ot
A^{\ot \bullet +1}, A^{\ot \bullet +1}\ot N);$ the homotopy
$h_{\bullet}$ on $B_{\bullet}^A(P)=P\ot A^{\ot \bullet +1}$ yields a
homotopy $- \circ h_{\bullet}$ from $id$ to $0$ on this
complex. Therefore, the double complex: 
$$\begin{array}{ccccc}
\vdots & & \vdots &\\
Hom_{A4}(B_1^A(P),A \ot N)
&\longrightarrow &Hom_{A4}(B_1^A(P),A \ot A \ot N)&\longrightarrow &\ldots  \\
\uparrow && \uparrow \\
Hom_{A4}(B_0^A(P),A \ot N) &\longrightarrow &Hom_{A4}(B_0^A(P),A \ot A \ot N)&\longrightarrow &\ldots  \\
\uparrow && \uparrow \\
Hom_{A4}(P,A \ot N) &\longrightarrow &Hom_{A4}(P,A \ot A \ot N)&\longrightarrow &\ldots  \\
\end{array}$$ which is the double complex $K^{\bullet  \bullet}(P,N)
$ to which we have added one row, has exact columns. Its homology is
therefore that of the first row $Hom_{A4}(P,A^{\ot \bullet +1} \ot N)
$ (\emph{cf.} \cite{W} p59-60). Since $P$ is projective and $A^{\ot
\bullet +1} \ot N$ is exact, the homology of this complex is zero in
positive degree.

\begin{lemma} the functor $H_{A4} ^* (-,N)$ takes short exact
sequences of Hopf bimodules $0 \rightarrow M' \rightarrow   M \rightarrow  M''
\rightarrow 0$ to long exact sequences of vector spaces $$\cdots \rightarrow H _{A4} ^n (M'',N) \rightarrow H _{A4} ^n (M,N)
\rightarrow H _{A4} ^n (M',N) \stackrel{\delta}{\rightarrow} H
_{A4} ^{n+1} (M'',N) {\rightarrow} \cdots $$\end{lemma}

Consider an exact sequence of Hopf bimodules $0 \rightarrow M' \rightarrow   M \rightarrow  M'' \rightarrow 0.$ Then for all $q \geq 0,$ the sequence $$(E_q) : 0 \rightarrow M'  \otimes A^{\otimes q} \rightarrow   M   \otimes A^{\otimes q} \rightarrow M'' \otimes A^{\otimes q} \rightarrow 0$$ is also an exact sequence of Hopf bimodules. 
Now any Hopf bimodule can be written $M \cong A \otimes M^R$ where $M^R = \{m \in M / \delta _R (m) = m \otimes 1 \}$ is the space of right coinvariants of $M,$ the left $A-$module and right $A-$comodule structures being given by $$\begin{array}{rcl}
b.(a \otimes v) & = & ba \otimes v \\
\delta _R \, (a \otimes v) & = & \sum a^{(1)} \otimes v \otimes a^{(2)},
\end{array}$$ and therefore any exact sequence of Hopf bimodules
splits as a sequence of left $A-$modules and right $A-$comodules,
since it yields an exact sequence of coinvariant spaces (a Hopf
bimodule is free and therefore injective as a right comodule, and $M^R
\cong Hom ^{-A}(k,M)$ where $k$ is the trivial right comodule). Thus,  the sequences $(E_q)$ split as sequences of left $A-$modules and right $A-$comodules. From this, we infer exact sequences \begin{eqnarray*}
0 \rightarrow Hom_{A-}^{-A} (M'' \otimes A^{\otimes q} , A^{\otimes p}  \otimes N) & \rightarrow &  Hom_{A-}^{-A} (M \otimes A^{\otimes q} , A^{\otimes p}  \otimes N) \\
 &  \rightarrow &  Hom_{A-}^{-A} (M' \otimes A^{\otimes q} , A^{\otimes p}  \otimes N) \rightarrow 0
\end{eqnarray*} for all non negative $p$ and $q.$ This leads to a long
exact sequence $$\ldots  \rightarrow H_{A4}^n (M',N) \rightarrow
H_{A4}^{n+1} (M'',N) \rightarrow H_{A4}^{n+1} (M ,N) \rightarrow
H_{A4}^{n+1} (M' ,N) \rightarrow \ldots \hspace{1cm} \square$$

\paragraph{Remark} In fact, we have proved that $H_{A4}^*(M,N)$ and
$Ext _X^* (M,N)$ are $\delta-$functors which are isomorphic, the
isomorphism being the unique map of $\delta-$functors extending the morphism $id_{Hom_{A4}(M,N)}.$ The same is true with the other variable.

\ 

We shall need, in theorem \ref{Yonedath}, the link between theorems \ref{Ospel} and \ref{Ospelth}:

\begin{prop}\label{degree1} The isomorphism of theorem \ref{Ospel} is equal to the isomorphism of theorem \ref{Ospelth} in degree 1.
\end{prop} 
\textit{Proof}: Let us call $\psi _N$ the isomorphism $H_{A4}^1(M,N) \tilde{\rightarrow} Ext _X^1 (M,N), $ and $\varphi _{MN}^*$ the isomorphism $H_{A4}^*(M,N) \tilde{\rightarrow} Ext _X^* (M,N). $ If we show that $\psi$ commutes with the first connecting homomorphism $\delta ^1,$ then there will be a unique map of  $\delta-$functors  $H_{A4}^*(M,-) \rightarrow Ext _X^* (M,-)$ extending the pair of maps ($id_{Hom_{A4}(M,-)}, \psi), $ which will necessarily be equal to $\varphi _M,$ since this is the unique map of $\delta-$functors extending $id_{Hom_{A4}(M,-)}.$ Hence $\psi$ will be equal to $\varphi _M^1.$

Therefore, let $(E) : 0 \rightarrow M' \stackrel{i}{\rightarrow} M' \oplus M'' \stackrel{p}{\rightarrow} M'' \rightarrow 0$ be an exact sequence of Hopf bimodules, where the Hopf bimodule structure of $M' \oplus M''$ is determined by a 1-cocycle $(\alpha , \beta  ) \in H^1_{A4} (M'',M').$ Let $g $ be a map in $Hom_{A4} (M,N).$ The connecting homomorphism in the long exact sequence associated to $(E)$ for the functor $Ext_X^* (M,-)$ is given by Yoneda multiplication on the right by $(E),$ so we want to prove that $\psi _{M''} (\delta ^1 g) = g \circ (E)$ (here $\circ$ is the Yoneda product). 

On the one hand, $g \circ (E)$ is the sequence of vector spaces $$0 \rightarrow N \rightarrow \; ^{\displaystyle{N \oplus M' \oplus M''}} / _{\displaystyle{<(g(m'),m',0),m' \in M'>}} \rightarrow M'' \rightarrow 0,$$ in which the Hopf bimodule stuctures of the middle space are given by 
$$(l,m',m'').a = (la, m'a + \alpha (m'' \otimes  a ) , m''a) =(la +g( \alpha ( m'' \otimes a )),m'a, m''a)$$ so that $(l,0,m'').a = (  la +g( \alpha ( m'' \otimes a )),0, m''a),$ and 
\begin{equation*}
\begin{split} \delta _L (l,m',m'') & = (\delta _L ^L (l) ,\delta _L ^{M'} (m') + \beta (m'') ,  \delta _L ^{M''}  (m'')) \\
& = (\delta _L ^L (l) + (1 \otimes g )  \beta (m''), \delta _L ^{M'}(m'),  \delta _L ^{M''}  (m''))
\end{split}
\end{equation*} so $\delta _L (l,0,m'') 
 = (\delta _L ^L (l) + (1 \otimes g )  \beta (m''),0,\delta _L ^{M''}  (m'')).$

Thus  $g \circ (E)$ is equivalent to $0 \rightarrow L \rightarrow L \oplus M'' 
\rightarrow  M'' \rightarrow 0, $ in which the Hopf bimodule structure is  determined by $(g  \alpha , (1 \otimes g )  \beta).$

On the other hand, $\delta ^1 g$ is a 1-cocycle $(k_0, k_1) \in H_{A4} ^1 (M'',N)$ such that $(k_0  (p \otimes 1), k_1  (p \otimes 1))= Dh,$ with $h \in Hom_{A-}^{-A} (M' \oplus M'' , N)$ satisfying $h  i =g.$ We want to show that $\delta ^1 g$ and $(g  \alpha , (1 \otimes g )  \beta)$ are equivalent 1-cocycles.

Let us write $$\begin{array}{rcl}
h:M' \oplus M'' & \longrightarrow & L \\
(m',m'') & \mapsto & g(m') +u(m'').
\end{array}$$

Then \begin{eqnarray*}
k_0(m'' \otimes a) &=& h((0,m'').a)-h(0,m'').a \\
 & =& h(\alpha (m'' \otimes a),m''a)-u(m'').a \\
 & =& g  \alpha (m'' \otimes a) + \delta_h u(m'' \otimes a).
\end{eqnarray*}
In the same fashion, we get $k_1 = (1 \otimes g)  \beta + \delta _c u. $ Therefore  $\delta ^1 g = (g  \alpha ,  (1 \otimes g ) \beta) + Du,$ which provides what we wanted. $\square$

\paragraph{Remark} The previous results give a unification of various cohomological theories of Hopf bimodules by viewing Hopf bimodules as $X-$modules. There is another way of looking at Hopf bimodules, that is as left modules over the Drinfel'd double $\mathcal{D}(A)$ of $A.$ Specifically, there is an equivalence of abelian categories $$\begin{array}{ccc}
X-\mathrm{mod} & \longleftrightarrow & \mathcal{D}(A)-\mathrm{mod} \\
M & \mapsto & M^R \\
A \otimes V & \leftarrow & V
\end{array}$$ which is exact ($ M^R$ is the space of right coinvariants). It follows that $$Ext_X^* (M,N) \cong  Ext_{\mathcal{D}(A)}^* (M^R,N^R).$$


\section{Cup-product on $H_{A4}^* (M,N)$}

Let $M$, $N$ and $L$ be Hopf bimodules over $A.$ There is a graded product $Ext_X^* (M,L) \otimes Ext_X^* (L,N) \rightarrow Ext_X^* (M,N),$ that is the Yoneda product.  Since the space $Ext_X^* (M,N)$ is isomorphic to $H_{A4}^* (M,N),$ this is also true for the Hopf bimodule cohomology. In this section, we shall establish a formula for a cup-product $$H_{A4}^* (M,L) \otimes H_{A4}^* (L,N) \rightarrow H_{A4}^* (M,N)$$ (which does not involve $X$), and prove that it corresponds to the Yoneda product.

\begin{prop} Let $f \in Hom_{A-}^{-A} (M \otimes A ^{\otimes p-s} ,A
^{\otimes s} \otimes L)$ be a $p-$cochain and $g \in Hom_{A-}^{-A} (L
\otimes A ^{\otimes q-r} ,A ^{\otimes r} \otimes N)$ be a
$q-$cochain. Set $n=p+q$ and $t=s+r.$ Define the $n-$cochain $f \smile
g \in Hom_{A-}^{-A} (M \otimes A ^{\otimes n-t} ,A ^{\otimes t}
\otimes N)$ by: \begin{align*} f \smallsmile g &(m \otimes \mathbf{a_{1, n-t}} )=\\
&\sum (-1)^{s(q-r)}(1 ^{\otimes s} \otimes g )\,[f(m \otimes
\mathbf{a_{1,p-s}} ) .(\Delta ^{(s-1)} (a^{(1)}_{p-s+1} \ldots a^{(1)}_{n-t} ) \otimes 1 ) 
 \otimes  \mathbf{a^{(2)}_{p-s+1,n-t}}]
\end{align*} where $\mathbf{a_{u,v}}=a_u \ot \ldots \ot a_v$ and
$\mathbf{a^{(2)}_{u,v}}=a^{(2)}_u \ot \ldots \ot a^{(2)}_v.$

The differential $D=d_h + d_c$ of the total complex associated to the
Hopf bimodule double complex is a right derivation for the cup-product
$\smile,$ that is $$D(f \smile g) = Df \smile g + (-1)^p f \smile
Dg,$$ so that the formula for $\smile$ yields a product $H_{A4}^* (M,L) \otimes H_{A4}^* (L,N) \rightarrow H_{A4}^* (M,N). $
\end{prop}
\textit{Proof}: Some straightforward computations, making use of the properties $\sum \Delta ^{(u)} (a^{(1)}) \otimes \Delta ^{(v)} (a^{(2)}) = \Delta ^{(u+v)} (a)$ and $\Delta _i \circ \Delta ^{(u)} = \Delta ^{(u+1)},$ lead to the following formulae:
$$\begin{array}{c}
d_h (f \smile g) = d_h f \smile g + (-1)^p f \smile d_h g \\ 
d_c (f \smile g) = d_c f \smile g + (-1)^p f \smile d_c g,
\end{array}$$ whence the result. $\square$ 

\paragraph{Remarks} 
(i) When $M$ and $N$ are equal to the Hopf bimodule $A,$ the cup-product can be written in the double complex defining $H_b^*(A,A):$ take $f $ in $Hom_k (A^{\otimes p-s} , A^{\otimes s})$ and $g$ in $ Hom_k (A^{\otimes q-r}, A^{\otimes r}).$ In this setting, their cup-product becomes \begin{eqnarray*}
(f \smile g) (a_1 \otimes \ldots \otimes a_{n-t}) & = & \sum  f(a^{(1)}_1 \otimes \ldots \otimes a^{(1)}_{p-s} ) \ \Delta ^{(s-1)} (a^{(1)}_{p-s+1} \ldots a^{(1)}_{n-t}) \\
& &\  \otimes \Delta ^{(r-1)} (a^{(2)}_{1} \ldots a^{(2)}_{p-s})\ g(a^{(2)}_{p-s+1} \otimes \ldots \otimes a^{(2)}_{n-t}). 
\end{eqnarray*}

(ii) Note that if $A=kG$ is a  group algebra, then $H_b^*(kG,kG)$ is the Hochschild cohomology $H\!H^*(kG,k),$ and the cup-product on $H_b^*(kG,kG)$ is the Hochschild cup-product on $H\!H^*(kG,k).$

\

Let us give a partial associativity property of this cup-product which will come in useful later:

\begin{prop}\label{assoc} Let $f, \  g$ and $h$ be cocycles
representing elements in  $H_{A4} ^p (M,L),$ $  H_{A4} ^q (L,N),$ and $ H_{A4} ^0 (M,N).$
Then $(f \smile g) \smile h = f \smile (g \smile h).$  
\end{prop}
\textit{Proof}: Indeed, both cocycles are elements of type $(u_t)_{t=0,\ldots , p+q}$ in $\bigoplus_{t=0}^{p+q} Hom_{A-}^{-A} (M \otimes A^{\otimes p+q-t}, A^{\otimes t} \otimes N),$ and  in both cases,  the degree $t$ component (that is $u_t$)  is equal to $$(f \smile g \smile h )_t = (1 ^{\otimes t} \otimes h) \circ (f \smile g)_t. \ \square$$

Let us now establish the link between this cup-product and the Yoneda product of extensions:

\begin{theo}\label{Yonedath} Let $M, \ N$ and $L$ be Hopf bimodules over $A.$ Let $\varphi _{MN} : H_{A4} ^* (M,N) \rightarrow Ext_X^* (M,N)$ be the isomorphism which extends $id_{Hom_{A4} (M,N)},$ and let $$\begin{array}{ccc}
Ext_X^p (M,L) \otimes Ext_X^q (L,N) & \longrightarrow & Ext_X^{p+q} (M,N) \\
F \otimes G & \mapsto & G \circ F 
\end{array} $$ be the Yoneda product of extensions.

Then, if $f$ is a $p-$cocycle and $g$ is a $q-$cocycle in the Hopf bimodule cohomology, the relationship between the products is given by $$\varphi _{MN}^{p+q} (f \smile g) = (-1)^{pq} \, \varphi _{LN}^q (g) \circ \varphi _{ML}^p (f):$$ the cup product and the Yoneda product are equal up to sign.  
\end{theo}

\paragraph{Remark} If the cup-product is replaced by $f \times g = (-1)^{\mid  f  \mid . \mid  g  \mid} f \smile g,$ the product $\times$ is equal to the Yoneda product, but $D$ is then a \textit{left} derivation for $\times .$

\

{\flushleft \textit{Proof of theorem \ref{Yonedath}}\,:}  Let $f$ be a
cocycle representing an element in $H_{A4}^p (M,L),$ set $F= \varphi
_{ML}^p (f),$ and consider the maps $\alpha _N^n =(-1)^{np} f \smile -
:H_{A4}^* (L,N) \rightarrow  H_{A4}^* (M,N) $ and the maps $\beta _N^n
= - \circ F : Ext_X^n (L,N)  \rightarrow Ext_X^{n+p} (M,N)$  for all
$n \in \mathbb{N}.$ We need to prove that the diagram \begin{equation}\label{diagram1}
\begin{CD} 
H_{A4}^n (L,N)@>{\alpha _N^n}>>H_{A4}^{n+p} (M,N) \\
     @V{\varphi _{-N}^n}VV                        @VV{\varphi _{-N}^{n+p}}V \\
Ext_X^n (L,N)@>{\beta _N^n}>>Ext_X^{n+p} (M,N)
\end{CD}
\end{equation} 
commutes for all $n \in \mathbb{N}.$

First, let us prove that $\alpha$ is a map of $\delta$-functors. Let
$0 \rightarrow N' \stackrel{u}{\rightarrow} N \stackrel{v}\rightarrow
N'' {\rightarrow}  0$ be an exact sequence of Hopf bimodules; it
induces the following diagram: $$\xymatrix@=8pt{
 &(L,N')_n \ar@{^{(}->}[rr]^-{(1 \ot u)\circ -}\ar@{.>}'[d][dd]^(-.5){D_n} & & (L,N)_n \ar@{->>}[rr]^-{(1 \ot v)\circ -}\ar@{.>}'[d][dd]^(-.5){D_n} & & (L,N'')_n  \ar@{->}[dd]^{D_n}\\
 (M,N')_{n+p} \ar@{<-}[ur]^{\alpha_{N'}^n} \ar@{^{(}->}[rr]^(.4){(1 \ot u)\circ -}\ar@{->}[dd]_{D_{n+p}} & & (M,N)_{n+p}  \ar@{<-}[ur]^{\alpha_{N}^n} \ar@{->>}[rr]^(.4){(1 \ot v)\circ -}\ar@{->}[dd]^(.3){D_{n+p}} & & (M,N'')_{n+p}
\ar@{<-}[ur]_(.6){\alpha_{N''}^n}\ar@{->}[dd]^(.3){D_{n+p}} \\
 & (L,N')_{n+1}\ar@{^{(}.>}'[r][rr]^(-.25){(1 \ot u)\circ -} & & (L,N)_{n+1} \ar@{.>>}'[r][rr]^(-.25){(1 \ot v)\circ -} & & (L,N'')_{n+1} \\ 
   (M,N')_{n+p+1} \ar@{^{(}->}[rr]^-{(1 \ot u)\circ -}\ar@{<.}[ur]^{\alpha_{N'}^{n+1}} & &  (M,N)_{n+p+1} \ar@{<.}[ur]^(.4){\alpha_{N}^{n+1}}
\ar@{->>}[rr]^-{(1 \ot v)\circ -} & & (M,N'')_{n+p+1}
   \ar@{<-}[ur]_(.6){\alpha_{N''}^{n+1}}\\ \\
}$$ where we have used the  notation $( L,N)_n
:=\bigoplus_{0\leq r \leq n} Hom _{A-}^{-A}(L\ot A^{\ot n-r},A^{\ot
r}\ot N).$

By construction of the connecting morphisms, to see that $\alpha$ is a
morphism of $\delta$-functors, it is enough to check that all the squares
in this diagram commute. The fact that the maps $(1 \ot u) \circ -$
and $(1 \ot v) \circ -$ commute with $\alpha$ is easily checked; they
also commute with the differential, using the fact that $u$ and $v$
are Hopf bimodule morphisms, and in particular that they are morphisms
of right modules and left comodules. Finally, $\alpha$ commutes with
the differential, since the latter is a derivation for the cup-product
and $f$ is a cocycle. Now the functor $H_{A4}^*(M,-)$
 is a universal $\delta$-functor, therefore $\alpha$ is the unique map
of $\delta$-functors extending $\alpha ^0.$

Now consider the maps $\beta _N^n.$ The functor $Ext_X^* (M,-)$ is a
universal $\delta-$functor, and $ \beta$ is a map of
$\delta-$functors (indeed, if $E$ is a  short exact sequence, the
connecting morphism in the cohomological exact sequence it induces for
the functor $Ext_X^* (L,-)$ is Yoneda multiplication on the right by
$E,$ see~\cite{BKI} Chapter 10 Proposition 5). Therefore, $\beta$ is
the unique map of $\delta-$functors  extending $\beta ^0; $ the map $  \beta ^0_N $ is the pushout of $F$ by $g$. 

Owing to the universal properties of the maps $\alpha$ and $\beta$, to
prove that the diagram (\ref{diagram1}) commutes, it is enough to
prove that the diagram
\begin{equation}\label{diag}
\begin{CD} 
H_{A4}^0 (L,N)@>{\alpha _N^0}>>H_{A4}^{p} (M,N) \\
       @V{id}VV                        @VV{\varphi _{-N}^p}V \\
Ext_X^0 (L,N)@>{\beta _N^0}>>Ext_X^{p} (M,N)
\end{CD}
\end{equation} commutes. We shall do it by induction on $p.$

Let us therefore first consider the case where $p=1.$ Then $f=(f_0 , f_1) \in Hom_{A-}^{-A} (M \otimes A ,L) \oplus Hom_{A-}^{-A} (M ,A \otimes L)$ is a 1-cocycle, and $F$ is the exact sequence $0 \rightarrow L  \rightarrow L \oplus M \rightarrow M \rightarrow 0$ in which the Hopf bimodule structures are determined by 
$$\begin{pmatrix} \mu _R ^L & f_0 \\ 0 & \mu _R ^M \end{pmatrix} \mathrm{\  and \  }
\begin{pmatrix} \delta _L ^L & f_1 \\ 0 & \delta _L ^M \end{pmatrix}.$$ 

If $g \in Hom_{A4} (M,L)$ is a 0-cocycle, then $f \smile g$ is a 1-cocycle, and $\varphi _{MN}^1 (f \smile g) $ is equal to $0 \rightarrow N \rightarrow N  \oplus M \rightarrow M \rightarrow 0$ where 

$$\begin{pmatrix} \mu _R ^N & (f \smallsmile g)_0 \\ 0 & \mu _R ^M \end{pmatrix} \mathrm{\  and \  }
\begin{pmatrix} \delta _L ^N & (f \smallsmile g)_1 \\ 0 & \delta _L ^M \end{pmatrix}$$ determine the Hopf bimodule structures. 

Moreover, as in the proof of proposition \ref{degree1}, $g \circ F$ is equivalent to the exact sequence $0 \rightarrow N \rightarrow N  \oplus M \rightarrow M \rightarrow 0$ in which the Hopf bimodule structures are determined by $g  f_0$ and $(1 \otimes g)  f_1.$ Since these are respectively equal to $(f \smile g)_0$ and $(f \smile g)_1,$ this proves that $g \circ F$ is equivalent to  $\varphi _{MN}^1 (f \smile g).$
This concludes the first case. It follows that for any $f \in H^1_{A4}(M,L)$ and any $g \in H^n_{A4}(L,N),$  $$ \varphi _{MN}^n (f \smile g) = (-1)^n  \varphi _{ML}^n (g) \circ F.$$

We can now study the general case: assume $f$ is a $(p+1)-$cocycle. Then we can write $F=E_p \circ E_1$ where $E_1$ is an extension in $Ext_X^1 (M,R)$ and $E_p$ is an extension in $Ext_X^p (R,L),$ for a Hopf bimodule $R.$
From the first case $(p=1),$ we have 
\begin{equation}\label{induct}
f= (\varphi_{MN}^{p+1})^{-1} (F) = (\varphi_{MN}^{p+1})^{-1} (E_p \circ E_1) = (-1)^p (\varphi_{MR}^1)^{-1} (E_1) \smile (\varphi_{RL}^p)^{-1} (E_p).
\end{equation}

Then, if $g$ is a 0-cocycle in $Hom_{A4}(L,N),$ we get

\begin{eqnarray*}
(\varphi_{MN}^{p+1})^{-1} (g \circ F) & = &  (\varphi_{MN}^{p+1})^{-1} (g \circ E_p \circ E_1) \\
&=& (-1)^p (\varphi_{MR}^1)^{-1} (E_1) \smile (\varphi_{RN}^p)^{-1} (g \circ E_p ) \\
&=& (-1)^p (\varphi_{MR}^1)^{-1} (E_1) \smile [(\varphi_{RL}^p)^{-1} ( E_p) \smile (\varphi_{LN}^0)^{-1} (g )] \\
&=& (-1)^p [(\varphi_{MR}^1)^{-1} (E_1) \smile (\varphi_{RL}^p)^{-1} ( E_p)] \smile (\varphi_{LN}^0)^{-1} (g ) \\
&=& f \smile g
\end{eqnarray*} using the first case, the induction hypothesis, proposition \ref{assoc} and relation (\ref{induct}).

Therefore the diagram (\ref{diag}) commutes in general, and this concludes the proof. $\square$

\end{document}